\documentstyle[12pt,euscript]{amsart}
\def\cal#1{{\EuScript#1}}

\newtheorem{defi}{Definition}
\newtheorem{lem}{Lemma}
\newtheorem{prop}{Proposition}

\newtheorem{cor}{Corollary}
\newtheorem{rem}{Remark}

\newcommand{\2}{\frac{1}{2}}

\def\GL{\mathrm{GL}}
\def\Ind{\mathop{\mathrm{Ind}}\nolimits}
\def\Frob{\mathop{\mathrm{Frob}}\nolimits}
\def\supp{\mathop{\mathrm{supp}}\nolimits}
\def\trace{\mathop{\mathrm{trace}}\nolimits}
\def\End{\mathop{\mathrm{End}}\nolimits}

\def\endproof{\unskip\nobreak\hfil\penalty50\hskip.5em\hbox{}\nobreak\hfill$\Box$}

{\makeatletter
\gdef\matrix#1{\null\,\vcenter{\normalbaselines\m@th
    \ialign{\hfil$##$\hfil&&\quad\hfil$##$\hfil\crcr
      \mathstrut\crcr\noalign{\kern-\baselineskip}
      #1\crcr\mathstrut\crcr\noalign{\kern-\baselineskip}}}\,}}

\begin{document}
\title[Harmonic Analysis on the Twisted  Finite Poincar\'e Half-Plane]
{Harmonic Analysis on the Twisted  Finite Poincar\'e Upper Half-Plane}

\author{Jorge Soto-Andrade}
\address{Jorge Soto-Andrade\\
Departamento de Matem\'aticas\\
Facultad de Ciencias\\
Universidad de Chile\\
Casilla 653\\
Santiago, Chile}
\email{sotoandr@@mate.uncor.edu,
sotoandr@@abello.dic.uchile.cl}

\author{Jorge Vargas}
\address{Jorge Vargas\\
FAMAF\\
Universidad Nacional de C\'ordoba\\
Ciudad Universitaria\\
C\'ordoba, Argentina}
\email{vargas@@mate.uncor.edu}

\thanks{Soto-Andrade was partially supported by FONDECYT
Grants 92-1041 and 1940590,   DTI -- U. Chile, ICTP, ECOS-France and NSF
Grant DMS-9022140 at MSRI. \ Vargas was
partially supported by
CONICET, CONICOR, SecytUNC, ICTP and TWAS}

\begin{abstract}
We prove that the induced representation from a non trivial character of
the Coxeter torus of $\GL (2,F) $, for a finite field $F$, is
multiplicity-free;  we give an explicit description of the corresponding
 (twisted) spherical functions and a version of the Heisenberg Uncertainty
 Principle.

\end{abstract}

\maketitle

\section{Introduction}
Let $F$ be a finite field, with $q$ elements, and $E $ be its unique
quadratic extension. Put $G = \GL (2,F)$ and denote by $K$ the Coxeter torus
of $G$, realized as the subgroup of all matrices $m_z   (z \in
E^\times)  $ of the
maps $ w \mapsto zw  (w \in E) $ with respect to a fixed $F$ - basis
of  $ E$. Recall that the finite homogeneous space $ {\cal H} = G/K$ may be
 looked upon
as the finite analogue of (the double cover of)  the classical Poincar\'e
Upper Half Plane (see \cite{arca}). Harmonic analysis on $ \cal
H$ amounts to decompose the induced representation $\Ind _K^G {\bf1} $ from
the unit character {\bf1} of $K$ to $G$.  We are interested here in the ``twisted''
version of this, i.e., the decomposition of the induced representation
$\Ind ^G_K \Phi$ from a non (necessarily) trivial character $\Phi$ of $K$ to
$G$. The real analogue of this case has been considered in \cite{eja}.
We prove that this representation is multiplicity-free, taking
advantage of the fact that this is so for $\Ind _K^G {\bf1}$ (see \cite{arca}) and
reducing the computation of the multiplicities in $\Ind ^G_K \Phi $ to the ones
in $\Ind _K^G {\bf1}$. We also give an explicit description of the corresponding
(twisted) spherical functions. Finally, we give a version of the
Heisenberg Uncertainty Principle

\section{The Multiplicity One Theorem for  $\Ind ^G_K \Phi$.}

\subsection {The case $\Phi = \bf 1$}
	We consider first the special case  $\Phi = \bf 1 $  in which the multiplicity
one theorem follows from a geometric argument. In fact, we have
$$
\Ind ^G_K \bf 1  \simeq  (L^2( \cal H), \tau),
$$
where $L^2(\cal H)) $ stands for the space of all complex functions
 on $\cal H$ endowed with the usual canonical scalar product, and  $\tau$
 denotes the natural representation of $G$ in
$L^2(\cal H) $, defined by  $(\tau_gf)(z) = f(g^{-1}.z) $,  where
$ z \mapsto  g.z  $  is  the homographic action of $G$ on $\cal H$,  given by
$$
g.z = \frac {az + b} {cz + d }   \hspace{3cm}
 \mbox{for}  g =  \left( \matrix{
a    &b\cr
c     &d \cr}\right) \in G,  z \in \cal H
.
$$
\begin{defi}
For all   $z,w \in \cal H $, we put
$$ D(z,w) = \frac{N(z - w)}{N(z - \bar w)} $$
with the convention that $D(z,w) = \infty $ if $w = \bar z$.
\end{defi}

\begin{prop}
$D$ is an orbit classifying invariant function for the homographic action of $G$ in $\cal H \times \cal H$.
\end{prop}
 \begin{cor}
The commuting algebra of $  (L^2( \cal H), \tau)$ is commutative.
\end{cor}
This follows from the fact that, the classifying invariant  $D$ being
 symmetric, the  $G$-orbits in  $\cal H \times \cal H$ are also symmetric.
\endproof

\subsection{The case of general $\Phi$}
Let's denote by $\phi$ the restriction of $\Phi$ to $F^\times$.
We will prove that every twisting  of an irreducible
representation  $  \pi^d_{\theta} $  of $G$ (where the superscript $d$ denotes the dimension
of $\pi$ and $\theta$ its character parameter) by the character $(\Phi +\Phi^q) $ is isomorphic to a
representation of the form $\pi^{d'}_{\theta'}+\pi^{d''}_{\theta''}$, when restricted to $K$. In
fact we will work with the characters $\chi^d_\theta$ of the irreducible
representations $\pi
^d_{\theta} $ of  $G$,
for which we keep the notations of \cite{tesis} or \cite{anna}.
\begin{lem}
On $K$ we have\\
 $(\Phi +\Phi^q) \chi^{q}_{\alpha, \alpha} = \chi^{q-1}_{\Phi (\alpha \circ N)} + \chi^{q+1}_{\phi \alpha, \alpha}$,\\
 $(\Phi +\Phi^q) \chi^{1}_{\alpha, \alpha} = \chi^{q+1}_{\phi \alpha, \alpha } - \chi^{q-1}_{\Phi (\alpha \circ N)}$,\\
 $(\Phi +\Phi^q) \chi^{q+1}_{\alpha, \beta} = \chi^{q+1}_{\phi \alpha,\beta} + \chi^{q+1}_{\alpha, \phi \beta}$,\\
 $(\Phi +\Phi^q) \chi^{q-1}_{\Lambda} = \chi^{q-1}_{\Phi \Lambda} +\chi^{q-1}_{\Phi^q \Lambda}.$
\endproof
\end{lem}

Now  for a character $\chi$ of $G$, we have
$ \chi \circ \Frob  = \chi$ on  $K$, as it follows from the character
table.
Therefore $ \sum_K \bar \Phi (k^q) \chi (k) = \sum_K \bar \Phi (k) \chi (k) $
because $\Frob $ is an involutive automorphism.

Hence, the multiplicity of $\pi$ in $\Ind _K^G \Phi $  equals
$ \2 \sum_K (\bar \Phi+\bar \Phi ^q) (k) \pi(k) $
and so it is just the average of the multiplicities in $\Ind ^G_K {\bf 1} $
of two representations
 of $G$ (one of
which  may be virtual!)
\begin{rem}
 Put $\pi^{q+1}_{\alpha, \alpha} =    \pi^q_{\alpha}  + \pi^1_{\alpha}$
and
$  \pi^{q-1}_{ \alpha \circ N} = \pi^q_{\alpha} -\pi^1_{\alpha}$
for every
$\alpha \in (F^\times)^{\wedge}$.
It is easy to check than in the degenerate cases $ \alpha = \beta$ (for
$\pi =
\pi^{q+1}_{\alpha, \beta}$ ) and
$\Lambda = \Lambda^q$ (for $\pi = \pi^{q-1}_\Lambda$) we find
 for the multiplicities $m_1 (\pi) $
\begin{equation}
m_1(\pi^{q+1}_{\alpha, \alpha}) = 1
 \hspace{4cm} (\alpha \in  (F^\times)^{\wedge})
\label{eq:mdegsp}
\end{equation}
and
\begin{equation}
m_1(\pi^{q-1}_{ \alpha \circ N}) = - \delta_{\alpha, 1}
  \hspace{4cm}(\alpha \in  (F^\times)^{\wedge} )
\label{eq:mdegsd}
\end{equation}
\end{rem}

Using the fact that the multiplicities of the irreducible representations
of $G$ in $\Ind _K^G \bf 1$ are at most one  and   also equations
 (\ref{eq:mdegsp}) and (\ref{eq:mdegsd}), we get that the multiplicities
 are also at most one in the more general    case of
$\Ind _K^G \Phi $.
\endproof

\subsection{The multiplicities  $ m_{\Theta,d}(\Phi)$  of $\pi^d_\Theta $
 in  $\Ind _K^G \Phi $
for general  $\Phi \in  (E^{\times})^{\wedge} $}

In  Table~1 below, $ \pi^d_\Theta$ denotes an irreducible representation of $G$, of dimension
  $d$ and parameter $\Theta$. Then $d \in \{ 1, q, q+1, q-1\}$ and   $\Theta$ is of the form
  $\{\alpha, \beta \} $, with  $ \alpha, \beta \in (F^{\times})^{\wedge}$ or
  $\{\Lambda, \Lambda^q \} $  with  $ \Lambda \in (E^{\times})^{\wedge}$
\begin{table}[hb]
\caption{The multiplicities   $ m_{\Theta,d}(\Phi)$}
\begin{center}
\begin{tabular}{|l|l|}
\hline
$\pi^d_\Theta $    &    $ m_{ \Theta, d }(\Phi)$   \\
\hline
$\pi^1_{\alpha, \alpha} $ & $ \delta_{\alpha^2, \phi} $ \\
\hline
$\pi^q_{\alpha, \alpha} $ & $\delta_{\alpha^2,\phi} -\delta_{\alpha \circ N,\Phi} $ \\
\hline
$\pi^{q+1}_{\alpha, \beta} $ & $\delta_{\alpha\beta,\phi} $   \\
\hline
$\pi^{q - 1}_{ \Lambda, \Lambda^q} $ & $\delta_{\lambda,\phi} - \delta_{\Lambda,\Phi} -
\delta_{\Lambda^q,\Phi}  $ \\
\hline
 \end{tabular}
 \end{center}
\end{table}

\noindent {\sc NOTATIONS.}
Here $ \alpha, \beta \in (F^{\times})^{\wedge}$   with  $ \alpha \neq \beta$  and
$ \Phi, \Lambda \in (E^{\times})^{\wedge}$  with   $ \Lambda \neq  \Lambda^q  $, and
$\lambda$ (resp. $\phi$ ) denotes the restriction of the  character $ \Lambda $ (resp.
$ \Phi$ ) to   $(F^{\times})^{\wedge}$.

\section{The twisted spherical functions}
\subsection{The averaging construction}
In this section $G$ denotes an arbitrary finite group,  $K$ a
 subgroup of $G$and $\Phi$ a one dimensional representation of $K$.
We notice that the spherical functions for the representation $\Ind _K^G \Phi $
are obtained as weighted  averages of the characters of  $G$. More precisely:
\begin{defi}
Let $ L^1(G)$  be the group algebra of  $G$, realized as the convolution algebra of all
complex functions of $G$ and let   $ L^1_\Phi(G,K)$  be the convolution algebra of  all complex
functions $f$ on  $G$ such that
$$  f(kgk') = \Phi(k)f(g)\Phi(k') $$
for all $g \in G,  k,k' \in K $.
For any $f \in L^1(G)$ put
$$
(P_ \Phi f)(g) = \frac{1}{|K|} \sum_{k \in K} \Phi^{-1}(k)f(kg)
$$
for all $g \in G$.
\end{defi}

Notice that the operator  $P_\Phi$ is just convolution with the idempotent function
$\varepsilon^\Phi_K \in L^1G  $ which coincides with  $|K|^{-1}\Phi$ on $K$ and vanishes
elsewhere. Moreover  $ L^1_\Phi(G,K) $ may be writen as  $
\varepsilon^\Phi_K \ast L^1G \ast \varepsilon^\Phi_K $ and its elements $f$ are
characterized by the properties
$$
      \varepsilon^\Phi_K \ast f = f = f \ast \varepsilon^\Phi_K
$$
\begin{lem}
Let $\chi$  be the character of an irreducible representation $\pi$ of
 $G$. Then
$P_{\Phi} (\chi) (e) \neq 0$  iff $\pi$ appears in $\Ind _K^G \Phi$.
\endproof
\end{lem}

\begin{lem}
$P_{\Phi} (\chi) $ is a non-zero function iff it doesn't vanish
for  $g=e$.
\endproof
\end{lem}

\begin{prop}
The mapping $P_ \Phi $ is an algebra epimorphism from the center $Z(L^1G))
$ of
the convolution
 algebra $L^1G $ onto the center $Z( L^1_\Phi(G,K))$
of the convolution  algebra   $ L^1_\Phi(G,K)$.
\end{prop}
{\em Proof:} We have
$$
\begin{array}{lllll}
P_ \Phi (f_1 \ast f_2) & = & \varepsilon^\Phi_K \ast (f_1 \ast f_2) & =
&( f_1 \ast \varepsilon^\Phi_K )\ast f_2   \\
& = & (f_1 \ast \varepsilon^\Phi_K \ast \varepsilon^\Phi_K )\ast f_2 &
=  & ( \varepsilon^\Phi_K \ast f_1 )\ast ( \varepsilon^\Phi_K  \ast f_2  ) \\
&  = & P_ \Phi f_1 \ast  P_ \Phi f_2.
\end{array}
$$
since $f_1$  is central and $\varepsilon^\Phi_K $ is idempotent. Moreover
the dimension $d$ of the image of $ Z(L^1G))$ under $P_ \Phi $ is the number
of irreducible characters $\chi$ of $G$ such that $P_\Phi(\chi) \neq    0$;
but   $P_\Phi (\chi) \neq   0$ iff  $(P_\Phi \chi)(e) \neq   0$ and, the
number $(P_\Phi (\chi))(e)$ being the multiplicity in $\Ind _K^G \Phi$
of the representation
$\pi$
of $G$ whose character is $\chi$, we see that $d$ is just the number
of irreducible representations $\pi$ of $G$ appearing in $\Ind _K^G \Phi$, i.
 e. the dimension of the center of $ L^1_\Phi(G,K)$.
\endproof

\begin{cor}
 The nonzero functions that satisfy the functional equation
 $$
 h(x)h(y) = \int_K \bar {\Phi }(k) h(xky) \,dk
 $$
 linearly span the center of the algebra $L^1_{\Phi} (G,K) $.

\end{cor}
 {\em Proof:}   The functions $h$ that satisfy the above functional equation are
 exactly the complex multiples of the functions $P_{\Phi} (\chi) $; for a
 proof (see \cite{4}). Therefore the corollary follows.
\endproof

\subsection{Explicit formulae for the twisted spherical functions}

Define
 $$
S^\Phi_\Lambda (a) = - (q^2 - 1)^{-1} \sum_{(z,w) \in
 \Gamma_a}\Phi^{-1}(z)\Lambda(w)
$$
 for  $ \Lambda \in  (E^\times)^{\wedge} $ and $a   \in
  F^\times$, where $\Gamma_a $ denotes
  the set of all $(z,w) \in
 E  ^{\times} \times E ^{\times}  $
 such that  $ N(w) = aN(z) $  and  $Tr(w) = 2(a + 1)^{-1} Tr(z) $.

Then the spherical function $ \zeta ^ \Phi_ \Lambda $ of $G$ associated to
 the cuspidal character $\chi^{q-1}_\Lambda $ of $G$ is given on the
 representatives  $d(a,1) = \left( \matrix{
 a   &   0 \cr
 0   &   1 \cr} \right) $ \hspace{0.5cm} $(a \in  F^\times)  $
 for the $K$ - double cosets in $G$,
 by
 \[
    \zeta ^ \Phi_ \Lambda (d(a,1))   =  S^\Phi_\Lambda (a)  +
   q(q + 1)^{-1} \delta_{a,1} \delta_{\lambda, \phi},
 \]
 where  $\lambda$ (resp. $\phi$) denotes the restriction of the character
  $\Lambda$ (resp. $ \Phi $) of  $E^ \times$  to  $ F^ \times$.

Notice that $a = 1$ corresponds to the origin in $\cal H$ and  $a = -1$
 corresponds to the antipode of the origin in $\cal H $.
 It is not difficult to check that these formulae for the spherical functions
 are equivalent to the ones given in \cite{arca} for the case  $\Phi = 1$.


\subsection{A new form for the cuspidal spherical functions for
  $\Phi = {\bf1 }$ ({\bf char $F$} $  \neq 2) $}

For $a \neq 1 $, one has the following new expression for the spherical functions
   estimated in \cite{katz}
  $$\zeta^\Phi_\Lambda (a) = (q + 1)^{-1} \sum_{u \in U} \varepsilon
  (Tr(u) - (a + a^{-1}))( \varepsilon \omega)(u),$$
 for  $a \neq 1$,  where  $\varepsilon$ denotes the
  sign character of   $ F^ \times$.

\section{Heisenberg Uncertainty Principle}
For this section, $G$ denotes an arbitrary finite group, $K$ any subgroup of $G$
and $\Phi $ any linear character of  $K$.

Let $ \hat G^{\Phi}$ be set of allthe equivalence classes of irreducible representations of $G$
that contain the character $\Phi$ when restricted to $K$. For each equivalence class we choose, once
and for all, a representative $(\pi,V_{\pi})$.
As usual, for each $f$ in $L^1(G)$, the Fourier Transform ${\cal F} (f) $, valued in the class
 $(\pi,V_{\pi})$, is the linear operator  $\cal F$  in $V_{\pi}$
defined by
$$
{\cal F} (f) (\pi):=\pi (f):=\frac{1}{|G|} \int_G f(g) \pi (g^{-1}) dg:=
\frac{1}{|G|}\sum_{ g \in G} f(g) \pi(g^{-1}). $$

We recall the statement of the Plancherel theorem for a function $ f \,\in L^1_{\Phi}(G,K) $
$$
f(g) =\frac{1}{G} \sum_{ \pi \in \hat G^{\Phi}} d_{\pi} \trace(\pi(f) \pi (g))  ; $$
here $ g\, \in \, G $ arbitrary and $d_{\pi}:=\dim V_{\pi} $.

For any complex valued function $f$ on $G$, let $ |\supp (f)| $ denote the number of elements of the support of $f$. That is, the number of points of $G$ where $f$ takes nonzero values.
\begin{prop}
[\bf Heisenberg Uncertainty Principle]
For any nonzero function $f \,\, \in \,\, L^1_{\Phi} (G,K) $ we have 
$$
|\supp (f) |\,( \sum_{\pi \in \supp ({\cal F} (f) )} d_{\pi})\, \geq\, |G|. $$
Here $\supp (\cal F (f) ) $ is the subset of $\hat G^{\Phi} $ where $\cal F (f)$ does not vanish.
\end{prop}
{\em Proof:} For any function $f$ on $G$ we recall that $$
\Vert f \Vert^2_{2} = \sum_{x \in G} \vert f(x) \vert ^2 ;\,\, \Vert f \Vert_{\infty} = \max_{ x \in G} \vert f(x) \vert ;\,\, \Vert f \Vert_2^2 \leq \Vert f \Vert^2_{\infty} \vert \supp (f) \vert \eqno (*)
$$

From now on, we fix a $G-$invariant inner product on $V_{\pi}$. Then $T^*$
denotes the adjoint of a linear operator $T$ on $V_{\pi}$ with respect to this inner product. Also
$\Vert T \Vert $ denotes the Hilbert-Schmidt norm on $\End  V_{\pi} $
 defined by $\trace(TS^*),$ for  $ S, T \in \End  V_{\pi} $

Since $ f \in L^1_{\Phi}(G,K) $, as we pointed out before,  
the Plancherel Theorem says that we have that 
$ \supp (f) $ is contained in $ \hat G^{\Phi} $ and that
$$
f(x) =\frac{1}{G} \sum_{ \pi \in \hat G^{\Phi}} d_{\pi} \trace(\pi(f)
\pi (x)).
$$

The Cauchy--Schwarz inequality applied to the Hilbert-Schmidt inner product says that the first of the
 two following inequallities is true,
$$
\trace (\pi (f) \pi (x) ) \leq \Vert \pi (f) \Vert \Vert \pi (x) \Vert \leq \Vert \pi (f) \Vert,
$$
the second inequality follows from the fact that $ \Vert T \Vert =1 $ for a unitary operator.

Putting together the last two statements we get $$
\Vert f \Vert_{\infty} \leq \frac{1}{G} \sum_{\pi \in \hat G^{\Phi}} d_{\pi} \Vert {\cal F} (f)(\pi) \Vert
$$
The classical Cauchy--Schwarz inequality and the fact
that $   d_{\pi} = d_{\pi}^{\frac 12}d_{\pi}^{\frac 12   }$  imply that
$$
\Vert f \Vert^2_{\infty} \leq \frac{1}{\vert G \vert^2} \sum_{\pi \in \hat G^{\Phi} } d_{\pi} \Vert {\cal F} (f)(\pi) \Vert^2 \,\, \sum_{\pi \in \supp({\cal F}(f))} d_{\pi}.
$$
Now the $L^2-$version of Plancherel Theorem says that
$$
\Vert f\Vert^2_2 =\frac{1}{\vert G \vert} \sum_{\pi \in \hat G^{\Phi}}
d_{\pi} \Vert {\cal F} (f)(\pi) \Vert^2. 
$$
Therefore,
$$
\Vert f\Vert_{\infty}^2 \leq \frac{1}{\vert G \vert} \Vert f \Vert^2 \,
\sum_{\pi \in \hat G^{\Phi}} d_{\pi} .
$$
Since $f$ is nonzero, we apply (*) to the above inequality and get the desired result.
\endproof

\end{document}